\documentclass[preprint,12pt]{elsarticle}
\usepackage{graphicx}
\usepackage{amssymb}
\usepackage{amsthm}
\usepackage{graphics}
\usepackage{epsfig}
\usepackage{amssymb}
\usepackage{hyperref}
\usepackage{graphicx}
\usepackage{subfigure}
\usepackage{float}
\usepackage{blkarray}
\usepackage{blkarray,bigstrut}
\usepackage{color}
\usepackage{tikz-network}
\usepackage{tikz}
\usepackage{tikz,pgfplots}
\pgfplotsset{compat=1.18}
\usetikzlibrary{patterns}
\usepackage[ruled,vlined]{algorithm2e}
\usepackage{caption}
\numberwithin{equation}{section}
\DeclareCaptionLabelFormat{cont}{#1~#2\alph{ContinuedFloat}}
\usepackage{fullpage}
\usepackage{amsmath}
\usepackage{amsfonts}
\usepackage{enumerate}
\usepackage{enumitem} 
\usepackage{tikz}
\usepackage{array}
\usepackage{bigdelim}
\usepackage{nicematrix}
\usepackage{mathtools}
\usepackage{blkarray, bigstrut}
\usepackage{gauss}
\usepackage{xcolor}
\usepackage{tabularray}
\usepackage{multicol}
\usepackage{multirow}
\usepackage{makecell}
\usepackage[most]{tcolorbox}
\UseTblrLibrary{booktabs}

\newcommand{\pf}{\noindent {\bf Proof: }}

\newtheorem*{theoremaux}{Theorem \theoremauxnum}
\gdef\theoremauxnum{1}

\newtheorem{lemma}{\bf Lemma}[section]
\newtheorem{theorem}{\bf Theorem}[section]

\theoremstyle{definition}

\journal{~}

\begin{document}

\begin{frontmatter}

\title{A characterization of a family of integral Cayley graphs}

\author[1]{Sauvik Poddar}
\ead{sauvikpoddar1997@gmail.com}

\author[2]{Sucharita Biswas\corref{cor1}}
\ead{biswas.sucharita56@gmail.com}


\address[1]{Department of Mathematics, Presidency University, Kolkata 700073, India\\
}
\address[2]{Department of Mathematics, Indian Institute of Technology Bombay, Mumbai 400076, India}
\cortext[cor1]{Corresponding author}

\begin{abstract}
Let $G$ be a finite group and $H$ be a subgroup of $G$. In our earlier work, we proved that the quasi-strongly regular Cayley graph $\Gamma_H(G)$ is integral whenever $H$ is normal in $G$, and raised the question of whether the converse also holds. In this paper, we answer this question affirmatively. We prove that if $\Gamma_H(G)$ is integral, then $H$ must be a normal subgroup of $G$. Consequently, $\Gamma_H(G)$ is integral if and only if $H\lhd G$. Our proof uses the irreducible representations of $G$ and the dimensions of their $H$-fixed subspaces. For a non-normal subgroup $H$, we construct a two-dimensional invariant subspace of a suitable representation operator and show that the corresponding eigenvalues are irrational.
\end{abstract}


\end{frontmatter}

\section{Introduction}

The study of the adjacency spectrum of a Cayley graph often reveals a close connection between the structure of a group and the properties of the graph. A graph is said to be \emph{integral} if all its adjacency eigenvalues are integers. The quest of characterizing integral graphs was initiated by Harary and Schwenk \cite{harary1974graphs}. Finding conditions under which a Cayley graph is integral is a well-studied problem in spectral graph theory. One can refer to \cite{alperin2012integral,cheng2019integral,cheng2023integral,godsil2025integral,huang2021integral,lu2018integral} in this context. In this paper, we consider a particular family of Cayley graphs arising from a finite group and one of its subgroups.

Let $G$ be a finite group with identity $1$ and let $H\le G$. Consider the set
$$S_H=\{(g,1),(1,g),(g,g)\mid g\in G\setminus H\}.$$
The graph $\Gamma_H(G)$ is the Cayley graph of $G\times G$ with respect to $S_H$, i.e.,
$$\Gamma_H(G)=\operatorname{Cay}(G\times G,S_H).$$
The set $S_H$ is inverse-closed, and hence $\Gamma_H(G)$ is an undirected graph. This family was introduced by Biswas and Das \cite{biswas2025family}, who studied its quasi-strong regularity, automorphism group, and transitivity properties. In particular, when $H$ is a non-trivial proper subgroup of $G$, the construction gives quasi-strongly regular graphs which are not strongly regular.

Recently, Poddar, Biswas, and Das \cite{poddar2026integral} studied the spectrum of $\Gamma_H(G)$ and proved that the graph is integral whenever $H$ is a normal subgroup of $G$. They also determined its complete spectrum in this case and showed that it depends only on $|G|$ and $[G:H]$. However, the question of whether normality of $H$ is also necessary for the integrality of $\Gamma_H(G)$ was left open \cite[Section 6]{poddar2026integral}. This question was motivated by computational examples in which non-normal subgroups gave rise to non-integral graphs.

In this paper, we solve this problem. Our main result is the following converse of the integrality theorem in \cite{poddar2026integral}.

\begin{theorem}\label{converse-theorem}
Let $G$ be a finite group and $H$ be a subgroup of $G$. Then the graph $\Gamma_H(G)$ is integral only if $H\lhd G$.
\end{theorem}

Combining Theorem \ref{converse-theorem} with the sufficient condition proved in \cite[Theorem 3.1]{poddar2026integral}, we obtain the complete characterization
$$\Gamma_H(G)\text{ is integral}\quad\Longleftrightarrow\quad H\lhd G.$$

Our proof is based on the representation theory of finite groups. We first show that if $H$ is not normal, then there exists an irreducible representation $\rho$ of $G$ whose $H$-fixed subspace is non-zero but is properly contained in the representation space. We then consider the operator associated with $S_H$ under the tensor product of $\rho$ and its dual representation. This operator has a two-dimensional invariant subspace, and we show that the two eigenvalues of the corresponding matrix are irrational. These are also eigenvalues of $\Gamma_H(G)$, proving that the graph cannot be integral.

The paper is organized as follows. In Section \ref{prelim-section}, we recall some results from linear algebra and representation theory needed for the proof. In Section \ref{necessary-section}, we establish the necessary condition for integrality and prove Theorem \ref{converse-theorem}.

\section{Preliminaries}\label{prelim-section}

Throughout the paper, all groups are finite and all vector spaces are finite-dimensional over $\mathbb{C}$. We use $H\le G$ to denote that $H$ is a subgroup of $G$, and $H\lhd G$ to denote that it is normal. We follow the notation for $\Gamma_H(G)$ used in \cite{poddar2026integral}.

\subsection{General preliminaries}

We begin with two elementary facts from linear algebra. The first allows us to express a tensor product operator as an operator on matrices. The second ensures that eigenvalues obtained by restricting an operator to an invariant subspace are eigenvalues of the original operator.

Let $M_{m\times n}(\mathbb{C})$ denote the set of all $m\times n$ matrices with complex entries. The \textit{vectorization} map $\operatorname{vec}:M_{m\times n}(\mathbb{C})\to\mathbb{C}^{mn}$ is a linear isomorphism defined as the vertical stacking of the columns of a matrix $A\in M_{m\times n}(\mathbb{C})$. Formally, for $A=[a_{ij}]$, we have 
$$\operatorname{vec}(A)=[a_{11},a_{21},\ldots,a_{m1},a_{12},a_{22},\ldots,a_{m2},\ldots,a_{1n},a_{2n},\ldots,a_{mn}]^T.$$

\begin{lemma}\cite[Lemma 4.3.1]{HornJohnson1991}\label{kronecker-identity-lemma}
Let $A\in M_{m\times n}(\mathbb{C}), B\in M_{p\times q}(\mathbb{C}), X\in M_{n\times p}(\mathbb{C})$. Then the following identity holds.
$$(B^T\otimes A)\operatorname{vec}(X)=\operatorname{vec}(AXB),$$
where $\otimes$ denotes the Kronecker product of matrices.
\end{lemma}

Let $V$ be a finite-dimensional vector space over $\mathbb{C}$ and $T:V\to V$ be a linear operator. An subspace $W$ of $V$ is said to be $T$-\textit{invariant} if $T(w)\in W$ for any $w\in W$.

\begin{lemma}\label{inv-subsp-eigenval-lemma}
Let $T:V\to V$ be a linear operator, and $W\subseteq V$ be a $T$-invariant subspace of $V$. Let $[T|_W]_{\mathcal{B}}$ be the matrix of the restricted operator $T|_W$ with respect to an ordered basis $\mathcal{B}$ of $W$. Then every eigenvalue of $[T|_W]_{\mathcal{B}}$ is also an eigenvalue of $T$ (and hence of the matrix representation $[T]_{\mathcal{B}'}$ with respect to any ordered basis $\mathcal{B}'$ of $V$).
\end{lemma}

\pf Let $\mathcal{B}=\lbrace w_1,w_2,\ldots,w_p\rbrace$ be an ordered basis of $W$. We extend $\mathcal{B}$ to an ordered basis $\mathcal{B}'=\lbrace w_1,w_2,\ldots,w_p,v_1,v_2,\ldots,v_q\rbrace$ of $V$. Since $W$ is $T$-invariant, $T(w_j)\in W$ for all $1\le j\le p$. Hence there exist scalars $c_{ij}\in\mathbb{C}$ ($1\le i,j\le p$) such that
\begin{align*}
T(w_1) &= c_{11}w_1+c_{21}w_2+\cdots+c_{p1}w_p+0\cdot v_1+0\cdot v_2+\cdots+0\cdot v_q,\\
T(w_2) &= c_{12}w_1+c_{22}w_2+\cdots+c_{p2}w_p+0\cdot v_1+0\cdot v_2+\cdots+0\cdot v_q,\\
&\quad \quad \quad \quad \quad \quad \quad \quad \quad \quad \quad\;\;\vdots\\
T(w_p) &= c_{1p}w_1+c_{2p}w_2+\cdots+c_{pp}w_p+0\cdot v_1+0\cdot v_2+\cdots+0\cdot v_q.
\end{align*}
Thus, the matrix representation of $T$ with respect to $\mathcal{B}'$ has the block upper-triangular form
$$[T]_{\mathcal{B}'}=\begin{pmatrix}
    [T|_W]_{\mathcal{B}} & B\\
    \mathbf{0} & A
\end{pmatrix},$$
where 
$$[T|_W]_{\mathcal{B}}=\begin{pmatrix}
    c_{11} & c_{12} & \cdots & c_{1p}\\
    c_{21} & c_{22} & \cdots & c_{2p}\\
    \vdots & \vdots & \ddots & \vdots\\
    c_{p1} & c_{p2} & \cdots & c_{pp}\\
\end{pmatrix},$$
$B$ and $A$ are matrices of sizes $p\times q$ and $q\times q$, respectively, and $\mathbf{0}$ is the $q\times p$ zero matrix. Denoting the characteristic polynomial of a square matrix $M$ by $\phi(M,x)$, we have
$$\phi([T]_{\mathcal{B}'},x)=\phi([T|_W]_{\mathcal{B}},x)\,\phi(A,x),$$
from which the result follows.\qed

\subsection{Representation and character preliminaries}

We now recall the basic representation-theoretic notations that will be used in the proof. We shall borrow the notations from \cite{poddar2026integral}. Let $G$ be a finite group and $V$ be a finite-dimensional vector space over $\mathbb{C}$. A \textit{representation} $(\rho,V)$ of $G$ is a group homomorphism $\rho:G\to GL(V)$, where $GL(V)$ denotes the group of all invertible linear transformations of $V$. We can simply denote a representation as $\rho$ if $V$ is already understood from the context. The \textit{degree} of $\rho$, denoted by $\deg(\rho)$, is $\dim(V)$, the dimension of $V$. For a representation $(\rho,V)$ of $G$, its character $\chi_\rho:G\to\mathbb{C}$ is defined as $\chi_{\rho}(g)=\operatorname{Tr}(\rho(g))$, where $\operatorname{Tr}(\rho(g))$ is the trace of the representation matrix of $\rho(g)$ with respect to some basis of $V$. We write $\operatorname{IRR}(G)$ for a complete set of inequivalent irreducible complex representations of $G$. For a subset $S\subseteq G$, we use the notation $\rho(S)\coloneqq\sum_{s\in S}\rho(s)$.
The group algebra $\mathbb{C}[G]$ is a vector space whose basis consists of the elements of $G$, i.e.,
$$\mathbb{C}[G]=\left\{\sum_{g\in G}c_gg\mid c_g\in\mathbb{C}\right\}.$$
The \textit{(left) regular representation} of $G$ is the homomorphism $L:G \to GL(\mathbb{C}[G])$ defined by
$$L(g)\left(\sum_{h\in G}c_hh\right)=\sum_{h\in G}c_hgh=\sum_{x\in G}c_{g^{-1}x}x,$$
for each $g\in G$. The character $\chi_L$ of the regular representation is given by
$$\chi_L(g)=\begin{cases}
    |G|,& \text{if $g=1$.}\\
    0, & \text{if $g\ne 1$.}
\end{cases}$$
Let $G$ be a group and $H\le G$. For $(\rho,V)\in\operatorname{IRR}(G)$, the $H$-fixed subspace of $V$ is defined as
$$V^H=\{v\in V\mid\rho(h)v=v,\forall h\in H\}.$$
whose dimension is given by
$$\dim(V^H)=\frac{1}{|H|}\sum_{h\in H}\chi_{\rho}(h),$$
which is equal to the rank (and trace) of the projection operator $P_H:V\to V$ where
$$P_H\coloneqq\frac{1}{|H|}\sum_{h\in H}\rho(h).$$
The image of the projection operator $P_H$ is $V^H$ and it satisfies $P_H^2=P_H$.

Let $(\rho,V)$ be a representation of $G$. The \textit{dual} $\rho^*$ of the representation $\rho$ is defined as,
$$\rho^*:G\to GL(V^*)$$
$$\rho^*(g)=\rho(g^{-1})^T,$$
where $V^*=\operatorname{Hom}(V,\mathbb{C})$ is the dual space of $V$. A representation is irreducible if and only if its dual is irreducible. Also for any representation $\rho$, $(\rho^*)^*=\rho$ and $\deg(\rho)=\deg(\rho^*)$. The \textit{kernel} of a representation $\rho$, defined as 
$$\ker(\rho)\coloneqq\lbrace{g\in G\mid \rho(g)=I}\rbrace$$
is a normal subgroup of $G$.

Let $G_1$ and $G_2$ be two groups and $(\rho_1,V_1)$, $(\rho_2,V_2)$ be two representations of $G_1$ and $G_2$ respectively. The \textit{tensor product} of the representations $(\rho_1,V_1)$ and $(\rho_2,V_2)$ is the representation $(\rho_1\otimes\rho_2,V_1\otimes V_2)$ of $G_1\times G_2$ defined as
$$\rho_1\otimes\rho_2:G_1\times G_2\to GL(V_1\otimes V_2)$$
$$(\rho_1\otimes\rho_2)(g_1,g_2)=\rho_1(g_1)\otimes\rho_2(g_2).$$
The irreducible representations of $G_1\times G_2$ are given by 
$$\operatorname{IRR}(G_1\times G_2)=\lbrace{\rho\otimes\psi\mid \rho\in\operatorname{IRR}(G_1),\psi\in\operatorname{IRR}(G_2)}\rbrace.$$
Let $\text{IRR}(G)=\{(\rho_i,V_i)\mid 1 \le i\le t\}$ be the complete set of inequivalent irreducible complex representations of $G$, where $\rho_1$ is the trivial representation. For each $1\le i\le t$, let $d_i=\deg(\rho_i)=\dim(V_i)$. The set of inequivalent irreducible representations of the group $G\times G$ is given by 
$$\operatorname{IRR}(G\times G)=\lbrace{\Omega_{ij}\coloneqq\rho_i\otimes\rho_j\mid 1\le i,j\le t}\rbrace.$$

\begin{lemma}\label{V^H=V-lemma}
Let $G$ be a group and $H\le G$. Let $(\rho,V)\in\operatorname{IRR}(G)$. Then $V^H=V$ if and only if $H\subseteq\ker(\rho)$.
\end{lemma}
\pf Suppose $V^H=V$. Let $h\in H$. Then 
$\rho(h)v=v$ for all $v\in V$, i.e., $\rho(h)=I$, which implies $h\in\ker(\rho)$. Hence $H\subseteq\ker(\rho)$.

Conversely, suppose $H\subseteq\ker(\rho)$. Let $v\in V$. Then $\rho(h)v=Iv=v$ for all $h\in H$. Hence $v\in V^H$. Thus $V\le V^H$. Since $V^H\le V$, we have $V^H=V$.\qed

\begin{lemma}\label{sum-rep-index-lemma}
Let $G$ be a group and $H\le G$. Then $$\sum_{i=2}^td_i\dim (V_i^H)=[G:H]-1.$$
\end{lemma}
\pf
\begin{align*}
\sum_{i=2}^td_i\dim (V_i^H)&=\sum_{i=2}^td_i\left(\frac{1}{|H|}\sum_{h\in H}\chi_i(h)\right)\\
&=\frac{1}{|H|}\sum_{h\in H}\left(\sum_{i=2}^td_i\chi_i(h)\right)\\
&=\frac{1}{|H|}\sum_{h\in H}\left(\sum_{i=1}^td_i\chi_i(h)-d_1\chi_1(h)\right)\\
&=\frac{1}{|H|}\sum_{h\in H}(\chi_L(h)-1),
\end{align*}
where $\chi_L$ is the character of the regular representation. Hence,
\begin{align*}
\sum_{i=2}^td_i\dim\left(V_i^H\right)&=\frac{1}{|H|}\left[(\chi_L(1)-1)+\sum_{h\in H,h\ne 1}(\chi_L(h)-1)\right]\\
&=\frac{1}{|H|}[(|G|-1)-(|H|-1)]\\
&=[G:H]-1.
\end{align*}
\qed

\medskip

\begin{lemma}\label{tensor-product-expression-lemma}
For any two non-trivial representations $\rho_i$ and $\rho_j$,
$$(\rho_i\otimes\rho_j)(S_H)=-|H|(P_i\otimes I_{d_j})-|H|(I_{d_i}\otimes P_j)+\sum_{g\in G\setminus H}\rho_i(g)\otimes\rho_j(g),$$
where $P_i=\frac{1}{|H|}\sum_{h\in H}\rho_i(h)$ is the projection operator from $V_i$ onto its $H$-fixed subspace $V_i^{H}$.
\end{lemma}

\pf Since $\rho_i$ and $\rho_j$ are non-trivial representations, by \cite[Lemma 4.4 (i)]{poddar2026integral}, $\rho_i(G)=\rho_j(G)=0$. From \cite[Eq. (3.3)]{poddar2026integral}, we have
\begin{align*}
(\rho_i\otimes\rho_j)(S_H)&=\sum_{g\in G\setminus H}(\rho_i\otimes\rho_j)(g,e)+\sum_{g\in G\setminus H}(\rho_i\otimes\rho_j)(e,g)+\sum_{g\in G\setminus H}(\rho_i\otimes\rho_j)(g,g)\\
&=\sum_{g\in G\setminus H}\rho_i(g)\otimes\rho_j(e)+\sum_{g\in G\setminus H}\rho_i(e)\otimes\rho_j(g)+\sum_{g\in G\setminus H}\rho_i(g)\otimes\rho_j(g)\\
&=\rho_i(G\setminus H)\otimes I_{d_j}+I_{d_i}\otimes\rho_j(G\setminus H)+\sum_{g\in G\setminus H}\rho_i(g)\otimes\rho_j(g)\\
&=(\rho_i(G)-\rho_i(H))\otimes I_{d_j}+I_{d_i}\otimes(\rho_j(G)-\rho_j(H))+\sum_{g\in G\setminus H}\rho_i(g)\otimes\rho_j(g)\\
&=(-|H|P_i)\otimes I_{d_j}+I_{d_i}\otimes(-|H|P_j)+\sum_{g\in G\setminus H}\rho_i(g)\otimes\rho_j(g)\\
&=-|H|(P_i\otimes I_{d_j})-|H|(I_{d_i}\otimes P_j)+\sum_{g\in G\setminus H}\rho_i(g)\otimes\rho_j(g).
\end{align*}
\qed

\section{Necessary condition for the integrality of the graph \texorpdfstring{$\Gamma_H(G)$}{GammaH(G)}}\label{necessary-section}

In this section, we prove that the integrality of $\Gamma_H(G)$ forces $H$ to be normal in $G$. The main step is to relate the normality of $H$ to the dimensions of the $H$-fixed subspaces of irreducible representations of $G$. We then use these dimensions to obtain an irrational adjacency eigenvalue when $H$ is not normal.

For convenience, we shall denote $|G|=n$, $|H|=k$, and $[G:H]=\ell=n/k$.  
For each $1\le i\le t$, let $r_i=\dim(V_i^H)$. We first prove a result which characterizes when $H\le G$ becomes a normal subgroup.

\begin{lemma}\label{H-normal-iff-characterization-lemma}
Let $G$ be a group and $H\le G$. Then $H\lhd G$ if and only if $r_i\in\{0, d_i\}$ for all $1\le i\le t$. Equivalently, if $H$ is not normal in $G$, there exists at least one $i\in\{2, \dots,t\}$ such that $0<r_i<d_i$.
\end{lemma}

\pf Let $H\lhd G$. Since $\dim(V_i^H)=\frac{1}{|H|}\sum_{h\in H}\chi_i(h)$, by \cite[Lemma 2.5 (ii)]{poddar2026integral}, we immediately have $r_i=d_i$, if $H\subseteq\ker(\rho_i)$ and $r_i=0$, otherwise.

Conversely, suppose $r_i\in\{0, d_i\}$ for all $1\le i\le t$. For each $i$ with $r_i=d_i$, we have $V_i^H=V_i$, which by Lemma \ref{V^H=V-lemma} implies $H\subseteq\ker(\rho_i)$. Define

$$N\coloneqq\bigcap_{\substack{1\le i\le t\\ r_i=d_i}}\ker(\rho_i).$$
Hence $N$ is a normal subgroup of $G$ containing $H$. Now 
\begin{itemize}
\item for every $i$ with $r_i=d_i$, we have $N\subseteq\ker(\rho_i)$, so by Lemma \ref{V^H=V-lemma} $\dim(V_i^N)=d_i=r_i$.
\item for every $i$ with $r_i=0$, since $H\le N$, we have $V_i^N\subseteq V_i^H =\{0\}$, so $\dim(V_i^N)=0=r_i$.
\end{itemize}
Thus $\dim(V_i^N)=r_i$ for all $1\le i\le t$. Now by Lemma \ref{sum-rep-index-lemma}, we get 
$$[G:H]-1=\sum_{i=2}^t d_i\dim(V_i^H)=\sum_{i=2}^t d_ir_i=\sum_{i=2}^t d_i\dim(V_i^N)=[G:N]-1,$$
which yields $|N|=|H|$. Since $H\le N$, we obtain $H=N$. This completes the proof.\qed 

\medskip


We now consider $H$ to be a non-normal subgroup of $G$. Then by Lemma \ref{H-normal-iff-characterization-lemma}, we can fix an index $i\in\{2,\dots,t\}$ such that $0<r_i<d_i$ and choose the corresponding irreducible representation $\rho_i$. Clearly, $\rho_i$ is a non-trivial and non-linear representation. To simplify notation, we write 
$$\rho\coloneqq\rho_i, \quad d\coloneqq d_i=\deg(\rho_i), \quad \text{and} \quad r\coloneqq r_i=\dim(V_i^H).$$
By Lemma \ref{kronecker-identity-lemma} and Lemma \ref{tensor-product-expression-lemma}, for any $X\in M_{d}(\mathbb
C)$, the action of the linear operator $\mathcal{M}\coloneqq(\rho^*\otimes\rho)(S_H):\mathbb{C}^{d^2}\to\mathbb{C}^{d^2}$ on $\operatorname{vec}(X)$ is given by
\begin{equation}\label{vec(X)-eqn}
\mathcal{M}(\operatorname{vec}(X))=-k\operatorname{vec}(XP)-k\operatorname{vec}(PX)+\sum_{g\in G\setminus H}\operatorname{vec}(\rho(g)X\rho(g^{-1})),
\end{equation}
where $P=\frac{1}{|H|}\sum_{h\in H}\rho(h)$ is the projection operator from $V_{\rho}$ onto its $H$-fixed subspace $V_{\rho}^{H}$.








\begin{lemma}
Let $v_I\coloneqq\operatorname{vec}(I_d)$ and $v_P\coloneqq\operatorname{vec}(P)$. Then the set $\lbrace{v_{I},v_P}\rbrace\subseteq\mathbb{C}^{d^2}$ is linearly independent.
\end{lemma}

\pf Since $\operatorname{vec}:M_d(\mathbb{C})\to\mathbb{C}^{d^2}$ is an isomorphism, showing the linear independence of $\lbrace{v_{I},v_P}\rbrace\subseteq\mathbb{C}^{d^2}$ is equivalent to showing the linear independence of $\lbrace{I_d,P}\rbrace\subseteq M_d(\mathbb{C})$. Let 
\begin{equation}\label{lin-indp-eqn}
c_1I_d+c_2P=0 \quad \text{for some} \quad c_1,c_2\in\mathbb{C}.
\end{equation}
Since $P^2=P$, we have from \eqref{lin-indp-eqn},
$$c_1I_d\cdot P+c_2P^2=(c_1+c_2)P=0.$$
Since $\operatorname{rank}(P)=r>0$, $P\neq 0$, which implies $c_1+c_2=0$, i.e., $c_2=-c_1$. Thus from \eqref{lin-indp-eqn},
$$c_1I_d-c_1P=c_1(I_d-P)=0.$$
As $\operatorname{rank}(P)=r<d=\operatorname{rank}(I_d)$, $P\neq I_d$. Hence $c_1=0$ which yields $c_2=0$.\qed

\begin{lemma}\label{mat-rep-M-lemma}
The $2$-dimensional subspace $W = \operatorname{span}\{v_I,v_P\} \subseteq\mathbb{C}^{d^2}$ is $\mathcal{M}$-invariant, and the matrix of the restricted operator $\mathcal{M}\vert{}_W$ with respect to the ordered basis $\{v_I,v_P\}$ is
$$M=\begin{pmatrix} 
(\ell-1)k & \frac{r}{d}\ell k\\ 
-2k & -3k 
\end{pmatrix}.$$  
\end{lemma}

\pf By Lemma \ref{lin-indp-eqn}, $W$ is a $2$-dimensional subspace of $\mathbb{C}^{d^2}$. To compute the matrix $M$, we compute $\mathcal{M}(v_I)$ and $\mathcal{M}(v_P)$, respectively, using \eqref{vec(X)-eqn}.
\begin{align*}
\mathcal{M}(v_I)
&=-k\operatorname{vec}(P)-k\operatorname{vec}(P)+\sum_{g\in G\setminus H}\operatorname{vec}(\rho(g)I_d\rho(g^{-1}))\\
&=-2k\operatorname{vec}(P)+\sum_{g\in G\setminus H}\operatorname{vec}(I_d)\\
&=((\ell-1)k)v_I+(-2k)v_P.
\end{align*}
Note that 
$$\sum_{g\in G \setminus H}\rho(g)P\rho(g)^{-1}=\sum_{g\in G}\rho(g)P\rho(g)^{-1}-\sum_{h\in H}\rho(h)P\rho(h)^{-1}=A-B.$$
The matrix $A$ is an intertwiner (i.e., it satisfies $\rho(x)A = A\rho(x)$ for all $x \in G$). Since $\rho$ is irreducible, Schur's lemma implies $A=cI_d$ for some scalar $c\in\mathbb{C}$. Taking traces yields
$$cd=\operatorname{Tr}(A)=\sum_{g\in G}\operatorname{Tr}(\rho(g)P\rho(g)^{-1})=n\operatorname{Tr}(P)=\ell k r \quad \text{i.e.,}\quad c=\frac{r}{d}\ell k.$$
For every $h \in H$, left and right multiplication by $h$ permutes $H$, so $\rho(h)P=P=P\rho(h)$. Hence $\rho(h)P\rho(h)^{-1}=P$ for all $h \in H$, giving $B=kP$.
\begin{align*}
\mathcal{M}(v_P)
&=-k\operatorname{vec}(P^2)-k\operatorname{vec}(P^2)+\sum_{g\in G\setminus H}\operatorname{vec}(\rho(g)P\rho(g^{-1}))\\
&=-2k\operatorname{vec}(P)+\operatorname{vec}\left(\frac{\ell kr}{d}I_d-kP\right)\\
&=\frac{\ell kr}{d}\operatorname{vec}(I_d)-3k\operatorname{vec}(P)\\
&=\left(\frac{r}{d}\ell k\right)v_I+(-3k)v_P.
\end{align*}
This establishes that $W$ is $\mathcal{M}$-invariant with representation matrix $M$.\qed

\medskip

Note that from Lemma \ref{sum-rep-index-lemma}, we have $\sum_{j=2}^{t}d_jr_j=\ell-1$. Also $d_jr_j\ge 0$ for all $1\le j\le t$. Hence we have
$$dr\le\ell-1 \quad \text{i.e.,} \quad \ell\ge dr+1.$$
Since $W$ is an $\mathcal{M}$-invariant subspace of $\mathbb{C}^{d^2}$, by Lemma \ref{inv-subsp-eigenval-lemma} the two eigenvalues of the matrix $M$ in Lemma \ref{mat-rep-M-lemma} are also eigenvalues of $\mathcal{M}$. By \cite[Eq. (3.2)]{poddar2026integral}, these eigenvalues are also the eigenvalues of the graph $\Gamma_H(G)$. We show that the eigenvalues of $M$ are always irrational whenever $0<r<d$ and $\ell\ge dr+1$.

\begin{lemma}\label{irrational-eval-M-lemma}
For $0<r<d$ and $\ell\ge dr+1$, the eigenvalues of $M$ are irrational.
\end{lemma}
\pf Suppose that $M$ has rational eigenvalues. The characteristic polynomial of $M$ is
$$\phi(M,x)=x^2-(\ell-4)kx+k^2\left(-3(\ell-1)+\frac{2r\ell}{d}\right),$$
whose discriminant is
$$\Delta=(\ell-4)^2k^2-4k^2\left(-3(\ell-1)+\frac{2r\ell}{d}\right)=k^2\left((\ell+2)^2-\frac{8r\ell}{d}\right).$$
Consider the integer 
$$D\coloneqq\frac{d^2}{k^2}\Delta=d^2(\ell+2)^2-8dr\ell.$$
Since $D=(d(\ell+2)-4r)^2+16r(d-r)$ and $0<r<d$, $D>0$. Also $\Delta$ is a square if and only if $D$ is a square. Note that $D$ is the discriminant of the polynomial
$$Q(y)\coloneqq y^2-d(\ell+2)y+2dr\ell.$$
Since $M$ has rational eigenvalues, the eigenvalues of $Q(y)$ are also rational. Evaluating $Q(y)$ at the two consecutive integers $y=2r$ and $y=2r-1$, we obtain,
$$Q(2r)=(2r)^2-d(\ell+2)(2r)+2dr\ell=4r^2-4dr=-4r(d-r)<0,$$
since $0<r<d$.
\begin{align*}
Q(2r-1)&=(2r-1)^2-d(\ell+2)(2r-1)+2dr\ell\\ 
&=4r^2-4r+1-4dr+d\ell+2d.
\end{align*}
Using the bound $\ell\ge dr + 1$, we obtain
\begin{align*} 
Q(2r-1)&\ge 4r^2-4r+1-4dr+d(dr+1)+2d\\ &=d^2r-4dr+4r^2+3d-4r+1\\ 
&=r(d-2)^2+4(r-1)^2+3(d-1)>0, \end{align*}
since $r\ge 1$ and $d\ge 2$. As $Q(2r)<0<Q(2r-1)$, there exists $y_0\in(2r-1,2r)$ such that $Q(y_0)=0$ with $y_0\in\mathbb{Q}$. Since $Q(y)$ is a monic polynomial with integer coefficients, by rational root theorem, $y_0\in\mathbb{Z}$, a contradiction. Hence $D\in\mathbb{Z}$ is not a perfect square and eventually $\sqrt{\Delta}=\frac{k}{d}\sqrt{D}\notin\mathbb{Q}$. Consequently, the eigenvalues of $M$ are irrational. This completes the proof.\qed

\medskip


\noindent\textbf{Proof of Theorem \ref{converse-theorem}:} The proof follows from Lemma \ref{H-normal-iff-characterization-lemma}, Lemma \ref{mat-rep-M-lemma} and Lemma \ref{irrational-eval-M-lemma}.\qed




\bibliographystyle{abbrv}
\bibliography{references}

\end{document}